\documentclass[11pt]{article}
\usepackage{amsmath,amsfonts, amssymb}
\usepackage[numbers]{natbib}
\usepackage{url}

\newcommand{\E}{\mathbf{E}}
\newcommand{\R}{\mathbb{R}}

\newcommand{\m}{\mathcal{M}}
\newcommand{\C}{\text{Cov}}

\newtheorem{lemma}{Lemma}
\newtheorem{thm}{Theorem}
\newtheorem{prop}{Proposition}
\newtheorem{cor}{Corollary}

  \newenvironment{Proof}{\noindent{\bf Proof} \ }{\QED}\smallskip
\newcommand\QED{\newline \rightline{$\blacksquare$} \bigskip}

 \smallskip
\newcommand\QEDa{\newline \rightline{$\blacksquare$} \smallskip}

\title{The Topology of Negatively Associated Distributions}
\author{Jonathan Root\footnote{ Department of Mathematics and Statistics, Boston University, Boston, MA 02215, USA}
, Mark Kon \footnote{ Department of Mathematics and Statistics, Boston University, Boston, MA 02215, USA. 
Research partially supported by NSF Award number DMS1736392}}
\date{}

\begin{document}

\maketitle

\begin{abstract}

We consider the sets of negatively associated (NA) and negatively correlated (NC) distributions as subsets of the
space $\mathcal{M}$ of all probability distributions on $\R^n$, in terms of their relative topological structures within the topological space of all measures on a given measurable space. We prove
that the class of NA distributions has a non-empty interior with respect to the topology of the total variation metric on $\mathcal{M}$.
We show however that this is not the case in the weak topology (i.e. the topology of convergence in distribution),
unless the underlying probability space is finite. We consider both the convexity and the connectedness of these classes of probability measures, and also consider the two classes on their (widely studied) restrictions to the Boolean cube in $\R^n$.
\end{abstract}

\section{Introduction}\label{sec:intro}

In recent years negatively associated probability distributions have been studied as potential generalizations of independent random variables \cite{NA,balls}.  However, the characterization of such probability measures on $\R^n$ has been elusive.  In many cases just the specialization of such a characterization to Boolean cube measures, i.e. probability measures whose marginals are simple variations of Bernoulli measures, has generated a great deal of interest \cite{pemantle,
ndgeometry}.  The characterization of the set of negatively associated measures on $\R^n$ can involve even simpler questions regarding the topological structure of this set within the space of all measures.  This question may have different answers under different topologies on the space of measures, which include the total variation topology and the standard weak (distributional) topology.  Simple versions of this question include whether the space of such distributions is connected, convex, closed, and whether it has an interior with respect to a given topology.


Denote by $\m(\R^n)$ the set of all Borel probability measures on $\R^n$. 
A probability measure $\mu\in \m(\R^n)$ is said to be {\it negatively correlated} (NC) if \begin{equation}\label{NC}
\int_{\R^n} x_ix_j \, d\mu(x) \leq \int_{\R^n} x_i \, d\mu(x) \int_{\R^n}x_j \,d\mu(x), \;\;\; \forall 1\leq i\neq j\leq n.
\end{equation}
We say that $\mu$ is {\it strictly} NC if strict inequality holds in (\ref{NC}).
We denote the class of NC distributions by $\m_{NC}$ or $\m_{NC}(\R^n)$. In this context, it is also common and equivalent to say that the variables themselves, $X_1,\dots, X_n$, are negatively correlated.

The functions $f_i(x)=x_i$, $i=1,\dots n$, are non-decreasing on $\R^n$. In general,
we say that a function $f:\R^n\to \R$ is non-decreasing
if $f(x) \geq f(y)$ whenever $x\geq y$ in the product ordering ($x\geq y$ if and only if $x_i\geq y_i$ for each $i=1,\dots, n$). We say that $f:\R^n\to \R$
is non-increasing if $f(x) \leq f(y)$ whenever $x\geq y$ (again in the product ordering). We
denote subsets of the index set $\{1,\dots, n\}$ by $I$, $J$, 
and define $x_I\in \R^{|I|}$ to be the restriction of $x\in \R^n$ to the index 
set $I$; here $|I|$  will denote the cardinality of $I$. Moreover, we denote by $\mu^{(I)}$ the marginal distributions
of $\mu$: for $A\subset \R^{|I|}$,
\begin{equation}\label{marginal}
\mu^{(I)}(A) := \int_A\int_{\R^{-I}} d\mu(x),
\end{equation}
where $\R^{-I}$ denotes all $x = (x_i)_{i\notin I}$.
Then $\mu$ is said to be {\it negatively associated} (NA) if for every 
disjoint $I,J \subset \{1,\dots, n\}$ and every non-decreasing and integrable
$f:\R^{|I|}\to \R$, $g: \R^{|J|}\to \R$, we have
\begin{equation}\label{NA}
\int_{\R^n} f(x_I)g(x_J) \, d\mu(x) \leq \int_{\R^n}f(x_I)\, d\mu(x) \int_{\R^n} 
g(x_J) \, d\mu(x),
\end{equation}
or equivalently
\begin{equation}\label{NA2}
\C_{\mu}(f(x_I),g(x_J)) \leq 0,
\end{equation}
where $\C_{\mu} :L^1(\R^n,\mu) \times L^1(\R^n,\mu) \to \R$ denotes
the covariance operator. Note that if $f$ or $g$ is constant, then we have trivial equality in (\ref{NA}). With that said, 
we say that $\mu$ is {\it strictly} NA if strict
inequality holds in (\ref{NA}) for all ($\mu$-almost surely) non-constant $f(x_I)$ and $g(x_J)$ and disjoint $I,J\subset\{1,\dots,n\}$. 
If we specify $f(x)=f_i(x)=x_i$ and $g(x)=g_j(x)=x_j$, $i\neq j$, then we arrive at (\ref{NC}). Thus negative association is stronger than negative correlation.
We denote the class of NA distributions
by $\m_{NA}$ or $\m_{NA}(\R^n).$ As in the case of negative correlation, it is common and equivalent to consider negatively associated variables. That is, if the variables $X_1,\dots, X_n$ are distributed according to a negatively associated distribution, then one may say that the variables themselves are negatively associated. 

Besides negative association, other attempts to quantify and conceptualize dependences among random variables appear in, for example, \cite{Block,Ebra,Jog,Karlin}.
A concept closely related to NA, known as positive association (PA), sheds some light on the class of NA distributions. The notion of positive association was introduced into the statistical literature prior to negative association, in \cite{Esary}. We say that $\mu$ is positively associated if $\C_{\mu}(f,g)\geq 0$, for all pairs of non-decreasing, real-valued functions $f$ and $g$. We note that we no longer assume that $f$ and $g$ are defined on disjoint subsets of variables, as we did with negative association. Remarkably (or not), significantly greater progress has been made in the theory of positive association than in the theory of negative association.  This, in part, is due to an elegant result known as the FKG inequality \cite{FKG}, which gives a sufficient criterion for PA. At its most basic level, the FKG inequality is known as Chebyshev's inequality \cite{Eaton} (distinct from the standard Chebyshev's inequality in elementary probability). This theorem states that if $X$ is a random variable on $\R$ (as opposed to $\R^n$), and $f,g:\R\to \R$ are both non-decreasing, then 
\begin{equation}\label{cheby}
\E (f(X)g(X)) \geq \E f(X)\E g(X).
\end{equation}
This holds for any probability distribution on the real line, so long as $f$ and $g$ are non-decreasing (or non-increasing). The proof of (\ref{cheby}) is straightforward, and follows from the basic pointwise inequality
\begin{equation}\label{cheby2}
(f(x)-f(y))(g(x)-g(y)) \geq 0,
\end{equation}
which holds for all non-decreasing (or non-increasing) $f,g :\R\to \R$. Indeed, assuming $x$ and $y$ are independent and identically distributed, upon expanding (\ref{cheby2}) and double integrating (over $x$ and $y$) we obtain Chebyshev's inequality (\ref{cheby}).

The FKG inequality is essentially a generalization of (\ref{cheby}) to the product setting, $\R^n$ equipped with the product ordering (i.e., $x=(x_1,\dots, x_n) \leq y=(y_1,\dots, y_n)$ iff $x_i\leq y_i$ $\forall i$). To state the result, we first define the functions $\wedge$ (meet, or greatest lower bound) and $\vee$ (join, or least upper bound) by,
\begin{align*}
x\wedge y &:= \max\{ z\in \R^n: z\leq x, z\leq y\}\\
x\vee y &:= \min\{ z\in \R^n : z\geq x, z\geq y\}.
\end{align*}
Then the FKG theorem states that if a discrete probability measure $\mu$ on $\R^n$ satisfies
\begin{equation}\label{fkg}
\mu(x\vee y)\mu(x\wedge y)\geq \mu(x)\mu(y)
\end{equation}
then $\mu$ is positively associated. 

Unfortunately, a criterion as simple as (\ref{fkg}) does not (yet) exists for negative association. As pointed out by Pemantle in \cite{pemantle}, the notion of negative association is not nearly as ``robust" as positive association. Since Chebyshev's inequality (\ref{cheby}) implies that any random variable is positively associated with itself, we cannot incorporate every non-decreasing function in the definition of negative association; but rather non-decreasing functions defined on disjoint coordinate subsets. And again, as noted in Pemantle \cite{pemantle}, there is a bound on how far $\E x_ix_j$ can lie below $\E x_i \E x_j$, due to the inequality $\text{Var}\left(\sum x_i\right) =\sum \C x_ix_j\geq 0$.

The study of the class of negatively associated random variables dates back to \cite{NA}, where basic properties of NA random variables are derived, and examples of NA random variables are given: multinomial, convolution of unlike multinomials, multivariate hypergeometric, Dirichlet, and Dirichlet compound multinomial variables. Though the notion of negative association has existed for some time, 
the literature on these distributions is still quite sparse \cite{NA,balls,pemantle,ndgeometry}. 
But interest in them is growing, due in part to the ease with which 
sums of NA (even NC) random variables satisfy sub-Gaussian tail bounds. Specifically, if $\mu$ is negatively correlated on $\R^n$, then 
\[
\mu\left( x\in \R^n: \left|\sum_{i=1}^nx_i - \E_{\mu}\sum_{i=1}^nx_i\right| \geq \lambda\right) \leq 
Ce^{-c\lambda^2},
\]
for some absolute constants $c,C>0$.
It is conjectured that the same may be true for the replacement of sums $\sum_i x_i$ by more general Lipschitz functions
of such variables; the only work on this question seems to come from \cite{peres}. (The notion that Lipschitz functions on a probability space concentrate about their mean, in the sense that their tails are (in the best case) sub-Gaussian, is called the concentration of measure phenomenon \cite{CM}.) The work in \cite{peres} seems to be inspired by the recent article \cite{ndgeometry}, in which the authors develop a novel notion of negative dependence known as the strong Rayleigh property. Their approach is via the geometry of associated generating polynomials, and they prove several conjectures put forth in this area of research.  

We emphasize here, however, that nowhere in the literature has the structure of the space of negatively associated or negatively correlated distributions been studied. It is therefore natural to ask about the topological or geometric properties of the space of NA, or even NC distributions. This question is the major impetus behind our work.

We consider these two classes of measures (NA and NC) broadly,
from a topological perspective. We view them as subsets of 
the general space of measures 
$\m(\R^n)=C_0(\R^n)^*$ (the dual of the space of continuous real-valued functions which vanish at infinity) endowed with the weak topology (technically this should be denoted as the weak-* topology).
This is the weakest topology ensuring the continuity of 
the maps $f\mapsto \int_{\R^n} f\, d\mu$ for $f\in C_0(\R^n)$, and coincides with the standard topology of convergence in distribution for measures. Thus we say that a
 sequence of distributions $\mu_n$ converges weakly to a distribution $\mu$
if,
\[
\int_{\R^n} f\, d\mu_n \to \int_{\R^n} f\, d\mu,
\]
for all $f\in C_0(\R^n)$. When $X\subset \R^n$ is compact in the standard topology, we may define
the weak topology on $\m(X)$ as follows. 
A basic open set in
the weak topology is given by \cite{parth, Reed}
\begin{equation}\label{openset}
V_{\mu}(f_1,\dots, f_k ; \epsilon_1,\dots, \epsilon_k) :=
\big\{\nu \in \m(X) : \left| \int f_i \, d\nu - \int f_i \, d\mu\right| < \epsilon_i, i=1\dots, k\big\}
\end{equation}
where $f_1,\dots, f_k$ are continuous real-valued functions on $X$. The family of sets obtained by varying 
$\mu, k$, $f_1,\dots, f_k$, $\epsilon_1,\dots, \epsilon_k$ form a basis for 
the weak topology, i.e. a collection of sets whose unions form all open sets. 
Thus a sequence of distributions $\mu_n$ converges weakly to 
a distribution $\mu$ if and only if 
\[
\int_{X} f\, d\mu_n \to \int_{X} f \, d\mu
\]
for every $f\in C(X)$ (now bounded due to compactness). 

We may in addition view the NC and NA families as subsets of the space of all measures $\m(\R^n)$, but now endowed with the total variation topology. This topology is induced from the  
{\it total variation distance}:
\[
\lVert \mu-\nu\rVert_{TV} := \sup_{|f|\leq 1} \left| \int_{\R^n} f
\, d\mu - \int_{\R^n} f \, d\nu\right|.
\]
In particular, in the setting of a discrete probability space (i.e. with support on a countable number of points), the total variation distance may be expressed as 
\begin{equation}\label{TVdiscrete}
\lVert \mu-\nu\rVert_{TV} = \sum_{x\in \R^n} |\mu\{x\} - \nu\{x\}|.
\end{equation}
Unless $\mu$ has finite support the total variation distance induces a stronger topology than the weak topology. This will be discussed later on in this paper.

An outline of the paper is as follows. 
We begin by showing
that the general class of NA distributions
on a compact subspace of $\R^n$ has a non-empty interior in the total variation topology, but not in the weak topology. We next specialize to 
the subspace of measures concentrated on $I_n=\{0,1\}^n$, the Boolean cube (a simplified space often considered [REFS?]), and consider
the interior and boundary of these distributions. This simple case affords intuitive arguments and constructive proofs. But it is still of great interest, and much is unknown about negative association on the Boolean cube \cite{pemantle,peres}.

Next we address the question of the convexity of the spaces of negatively associated and negatively correlated distributions. We show that these spaces are not convex for distributions on $\R^n$, and they are similarly non-convex when restricted to the Boolean cube. We then address whether or not these spaces are connected in the weak or total variation topology.

\section{The Topology of $\m_{NC}$ and $\m_{NA}$}

\subsection{The Interior of $\m_{NC}(\R^n)$ and $\m_{NA}(\R^n)$}

The interiors of both the space of negatively correlated and negatively associated distributions are intimately connected to their {\it strict} counterparts.  Recall that a distribution $\mu$ on $\R^n$ is strictly NC if
\[
\C_{\mu}(x_i,x_j) <0
\]
for all $1\leq i< j \leq n$, and a distribution $\mu$ on $\R^n$ is {\it strictly} NA if 
\[
\C_{\mu}(f(x_I),g(x_J)) <0 
\]
for all strictly non-decreasing (not almost surely constant) $f,g$ and disjoint $I,J \subset\{1,\dots, n\}$. If a distribution $\mu$ is strictly NA then it must be strictly NC. It is not {\it a priori} clear that strictly NA distributions even exist. However, we prove below in Lemma \ref{strictNAexist} that they indeed exist on the Boolean cube $\{0,1\}^n$ (and thus by extension on $\R^n$).

We note that here and elsewhere, the notion of a strictly non-decreasing function $f$ means by implication that $f$ is {\it strictly} non-decreasing, i.e. that it is not essentially constant with respect to the measure $\mu$ under consideration. We define the total variation of a non-decreasing function $f$ with respect to a measure $\mu$ to be $\sup f -\inf f$, where sup denotes essential sup (i.e. modulo sets of measure 0) and inf denotes essential inf. Note however that on the cube these two notions (sup and essential sup) coincide for non-decreasing and non-increasing functions, as do inf and essential inf.

\begin{lemma}\label{strictNAexist}
There exist strictly NA distributions on $\{0,1\}^n$
such that for all non-decreasing $f(x_I),g(x_J)$ (with $I,J$ disjoint sets of indices) having total variation 1 on $\{0,1\}^n$, there is an $\epsilon>0$ such that
\[
\int fg d\mu \leq \int f d\mu \int g d\mu  -\epsilon .
\]
Consequently, there exist strictly NA distributions on the whole of $\R^n$ satisfying the above equation (under the same measure $\mu$ supported on $\{0,1\}^n$ viewed as a subset of $\R^n$).
\end{lemma}
\begin{Proof}
Let $I_{n,1}$ denote the collection of vectors $(x_1,\dots, x_n) \in \{0,1\}^n$ such that $\sum_i x_i = 1$, and let $\mu$ be any probability distribution supported on $I_{n,1}$. Thus for some $\epsilon>0$, we have $\mu(x)>\sqrt{\epsilon}$ for all $x\in I_{n,1}$. 
Now assume that, as stated, $f(x_I)$ and $g(x_J)$ are non-decreasing functions of total variation 1, with $I$ and $J$ disjoint subsets of ${1,\ldots, n}$. First note that to check the condition 
\[
\C_{\mu}(f(x_I),g(x_J)) <0,
\]
it suffices to assume that $f(0,\dots,0) = 0$ and $g(0,\dots, 0) = 0$. Indeed, we may replace $f(x_I)$ with $f(x_I) - f(0,\dots, 0)$ and $g(x_J)$ with $g(x_J) - g(0,\dots, 0)$ without changing $\C_{\mu} (f(x_I),g(x_J))$. 
It follows that 
\[
\E_{\mu}[ f(x_I)g(x_J)]=0,
\]
since for $x$ in the support of $\mu$, $\sum_ix_i=1$ and so we must have $x_I=(0,\dots, 0)$ or
$x_J=(0,\dots, 0)$. On the other hand, since $f$ and $g$ are strictly non-decreasing (i.e. non-constant), zero at the zero vector, and of total variation 1, each must equal one at one or more points in the support of our measure, i.e., in $I_{n,1}$.
Thus $\E_{\mu}[f(x_I)]>\sqrt{\epsilon}$ and $\E_{\mu}[g(x_J)] >\sqrt{\epsilon}$, and so 
\[
\E_{\mu}[ f(x_I)g(x_J)]=0 <
(\sqrt{\epsilon})^2=\epsilon\leq
\E_{\mu}[f(x_I)]\E_{\mu}[g(x_J)].
\]
Therefore $\mu$ is strictly negatively associated satisfying the lower $\epsilon$-bound in the statement of the theorem.
We note that since $\mu$ satisfies this bound as a measure on the cube $I_n$, it also satisfies this bound when viewed as a measure on $\R^n$ (that is concentrated on $I_n\subset \R^n$).
\end{Proof}

We now move to the main results of this section. We begin by considering
the weak interior (i.e. interior in the weak topology on measures) of the space of negatively associated distributions on $\R^n$ or on any fixed open subset $G\subset \R^n$. The following Proposition shows that the weak interior of the NA distributions on $G$ is in fact empty.

\begin{prop}\label{empty}
Consider the space of probability distributions supported on a fixed open set $G\subset \R^n$ (or all of $\R^n$).
If $\mu$ is strictly NA on $G$, then every weak neighborhood of $\mu$ contains a non-NA distribution $\nu$.
\end{prop}

\begin{Proof}
Let $\mu$ be strictly NA and let $\epsilon$ and continuous bounded $f_1,\dots, f_k$ be given. Then $\nu'$
will be in the weak neighborhood $V_\mu(\epsilon;f_1,\dots, f_k)$
if and only if 
\[
\left| \int f_i \, d\mu - \int f_i \, d\nu'\right| \leq \epsilon
\]
for every $i=1\dots, k$. 
We construct a non-negatively associated distribution $\nu'$ in the weak neighborhood $V_{\mu}(\epsilon;f_1,\dots,f_k)$. This is done by way of a discrete distribution $\nu'$ satisfying $\int f_i d\mu =\int f_i d\nu'$ for each $i=1\dots k$, but which itself is not negatively associated. 

Let $\nu$ be a positively associated distribution on $G$
such that 
\[
\int f_{k+1}f_{k+2} \, d\nu  - \int f_{k+1}\, d\nu \int f_{k+2} \, d\nu >0,
\]
for some non-decreasing $f_{k+1},f_{k+2}$, and append the latter two functions to the above sequence, yielding $f_1,\ldots,f_{k+2}$. (See \cite{FKG} for existence theorems and examples of positively associated measures.)
Assume without loss of generality that $f_1,\dots,f_{k+1}, f_{k+2}$, $f_{k+1}f_{k+2}$
are linearly independent. For each $x$ consider the vector 
\[
{\bf f}_x=(f_1(x),\dots, f_k(x),f_{k+1}(x),f_{k+2}(x), f_{k+1}(x)f_{k+2}(x)) \in \R^{k+3},
\] 
with the last entry a product of $f_{k+1}$ and $f_{k+2}$. The collection $\{{\bf f}_x\}_{x \in G}$ spans $\R^{k+3}$. 
Thus we can find $\alpha_1,\dots, \alpha_{k+3}$ and $x_1,\dots, x_{k+3}$
so that 
\[
(\mu(f_1),\dots,\mu(f_k), \nu(f_{k+1}), \nu(f_{k+2}),\nu(f_{k+1}f_{k+2})) = \sum_{j=1}^{k+3}\alpha_j{\bf f}_{x_j}
\]
where $\mu(f_i) := \int f_i \, d\mu$. 
Therefore 
\[
\mu(f_i) = \sum_{j=1}^{k+3}\alpha_jf_i(x_j)
\]
 for each $i=1,\dots, k$,
 \[
 \nu(f_{k+1}) = \sum_{j=1}^{k+3} \alpha_j f_{k+1}(x_j),
 \]
  \[
 \nu(f_{k+2}) = \sum_{j=1}^{k+3} \alpha_j f_{k+2}(x_j),
 \]
 and
\[
\nu(f_{k+1}f_{k+2}) = \sum_{j=1}^{k+3} \alpha_j f_{k+1}(x_j)f_{k+2}(x_j).
\]
So the discrete distribution $\nu' = \sum_{j=1}^{k+3} \alpha_j \delta_{x_j}$
is in the weak neighborhood $V_{\mu}(\epsilon ; f_1,\dots, f_k)$, but it is
not negatively associated.
\end{Proof}

The study of the interior of the set of negatively associated distributions is complicated by the fact 
that the covariance condition for negative association must be checked on infinitely many
functions (in order to establish negative association for a single measure). This situation can however be avoided when the distribution is supported on a finite
subset of $\R^n$
(by a compactness argument in section \ref{sec_interior} below).
On a finite product probability space $X^n$, the space of probability
distributions is finite dimensional, and we may conclude as will be done in 
section \ref{sec_interior} that the interior of the collection of NA distributions is non-empty. Note that since the set of distributions on a finite space is finite dimensional, the two topologies (weak and total variation) discussed here coincide.

On the other hand, a probability measure is negatively correlated (on $\R^n$) if finitely many covariance conditions are satisfied (\ref{NC}). Because of this, we may prove that the weak interior of the class of NC distributions is non-empty in the space of distributions supported on a fixed compact set $X\subset \R^n$.


\begin{prop}
Let $X\subset \R^n$ be a compact subset. Then  $\m_{NC}(X)\subset \m(\R^n)$ has a non-empty interior in the weak topology.   
\end{prop}

\begin{Proof}
It suffices to show that $\int x_ix_j d\mu-\int x_i d\mu\int x_jd\mu$ is continuous in $\mu$, which follows from the fact that the functions $f_{ij}(x) = x_ix_j$, $1\leq i,j\leq n$,
and $f_k(x) = x_k$, $1\leq k\leq n$ are bounded.
 \end{Proof}
 
 This result however does not hold when we consider the class of negatively correlated distributions supported on all of $\R^n$. 

\begin{prop}\label{NCcompactinterior}
The collection of negatively correlated distributions on $\R^n$ has no interior in the total variation topology (and hence in the weak topology).
\end{prop}
\begin{Proof}
Let $\mu$ be a negatively correlated distribution. Consider the distribution $\nu_c=\frac{1}{2}(\delta_{-c{\bf 1}} +\delta_{c{\bf 1}})$ (a sum of point masses at two points), where $c>0$ is large and ${\bf 1} =(1,1,\dots,1)$. We note that $\nu_c$ has total variation 1. We claim that for any neighborhood $V$ of $\mu$ (in the TV metric), there is a distribution in $V$ that is positively correlated, of the form $\mu_{\alpha,c}= \alpha\mu +(1-\alpha)\nu_c$, for some $\alpha\in (0,1)$ close to 1 and $c>0$.

The idea here is that the distribution $\nu_c$ has a positive correlation that is arbitrarily large as $c$ becomes large, so that adding only a small multiple $(1-\alpha)v_c$ (if $c$ is large) will cause a distribution to become positively correlated. 

Note first that
\[
\lVert \mu - \mu_{\alpha,c}\rVert_{TV} = \lVert (1-\alpha)\mu +(\alpha-1)\nu_c\rVert_{TV}
\]
\[
\leq (1-\alpha)\lVert \mu \rVert_{TV} +(1-\alpha)\lVert \nu_c \rVert_{TV} = 2(1-\alpha)
\]
which is arbitrarily small for $\alpha$ close to 1 (uniformly in $c$). Hence, uniformly in $c$, the measure $\mu_{\alpha,c}$ is in $V$ for $\alpha$ sufficiently close to 1, which we assume is the case.
However for this value of $\alpha$, we now allow $c$ to grow larger. Note that the covariance
\[
\int x_i x_j d\mu_{\alpha,c}
= \int x_ix_jd[\alpha\mu +(1-\alpha) \nu_c]
\]
\[
    =\alpha\int x_ix_j d\mu +\frac{1}{2}(1-\alpha) \int x_ix_j d[\delta_{-c{\bf 1}} +\delta_{c{\bf 1}}]
    \]
    \[
    =\alpha\int x_ix_j d\mu +\frac{1}{2}(1-\alpha)[c^2+c^2],
    \]
which for sufficiently large $c$ is clearly positive. Thus any TV neighborhood $V$ of $\mu$ has a positively correlated distribution in it.
\end{Proof}






As a consequence of Proposition \ref{NCcompactinterior} we have

\begin{prop}
The space $\m_{NA}$ on all of $\R^n$ has an empty interior with respect to both the TV and weak topologies.
\end{prop}

However, on a compact set $X\subset \R^n$, the space $\m_{NA}(X)$ has a non-empty TV interior, as is shown here:

\begin{prop}
Let $X\subset \R^n$ be a compact subset. Then
$\m_{NA}(X)$ has a non-empty interior with respect to the 
total variation metric.
\end{prop}

\begin{Proof}
It is not hard to show that if $\lVert \mu - \nu\rVert_{TV} <\epsilon$,
then $|\C_{\mu}(f,g) - \C_{\nu}(f,g) | \leq 3\epsilon$ for every $f,g$ satisfying
$\lVert f\rVert_{\infty},\lVert g\rVert_{\infty} \leq 1$.

According to Lemma \ref{strictNAexist}, let $\mu$ be strictly NA so that $\C_{\mu}(f,g) < -\epsilon$ for every strictly non-decreasing $f(x_I),g(x_J)$ of total variation 1 defined on disjoint index subsets $I,J\subset\{1,\dots, n\}$, and some $\epsilon>0$.
Choose
$\mu'$ such that $\| \mu - \mu' \|_{TV}<\epsilon/6$. Then
\begin{equation}\label{TV}
\C_{\mu'}(f,g) < \epsilon/2 + \C_{\mu}(f,g) < -\epsilon/2 <0.
\end{equation}
For general $f$ and $g$ (not of total variation 1) (\ref{TV}) holds by multiplying these by constants, without changing the negative covariance. This completes the proof. 
\end{Proof}

For a compact $X$, since $\m_{NA}(X) \subset \m_{NC}(X)$, it immediately follows that the interior of the collection of NC distributions is non-empty in the total variation topology (this also follows from the fact that it is non-empty in the weak topology):
\begin{cor}
For a compact $X\subset \R^n$,
$\m_{NC}(X)$ has a non-empty interior with respect to the total variation metric.
\end{cor}

\subsection{$\m_{NC}$ and $\m_{NA}$ on the Boolean Cube}\label{sec_interior}


We reformulate the conditions for negative correlation and negative association on the Boolean cube $I_n=\{0,1\}^n$ as polynomial inequalities. We consider a compactness argument in the case of strict negative association in order to get a handle on the infinity of conditions contained within  definition (\ref{NA}). Restricting these measures to the Boolean cube affords great flexibility due to the topological properties of both the space of probability measures and the space of continuous functions on said cube.

Denote by $\mu^{(i)}$, $i=1,\dots,n$,
and $\mu^{(i,j)}$, $1\leq i,j\leq n$,
respectively, the one and two-dimensional marginals of $\mu$. 
On the Boolean cube $I_n=\{0,1\}^n$ we have $\E x_i =\mu^{(i)}(1)$, for each $i=1,\dots, n$, and $\E x_ix_j = \mu^{(i,j)}(1,1)$, for each $1\leq i,j\leq n$. The condition for negative correlation therefore reduces to,
\begin{equation}\label{NCBoolean}
\mu^{(i,j)}(1,1)\leq \mu^{(i)}(1)\mu^{(j)}(1), \;\;\; \forall 1\leq i,j\leq n.
\end{equation}
Any probability measure $\mu$ on the Boolean cube is uniquely determined
by a vector of length $2^n$, $\mu = (\mu_1,\mu_2,\dots, \mu_{2^n})$, 
such that $\sum_i \mu_i =1$. Thus equation (\ref{NCBoolean}) 
may be written in the form 
\begin{equation}\label{NCPolynomial}
p_{\mu}(\mu_1,\dots,\mu_{2^n}) \leq 0,
\end{equation}
where $p_{\mu}$ is a polynomial in $\mu_1,\dots,\mu_{2^n}$.

Now consider definition (\ref{NA}) of negative association, specialized to the Boolean cube. Denote by $I$ and $J$ disjoint subsets of indices in $\{1,\dots, n\}$,
and by $x_I,x_J$ vectors restricted to the indices of $I$ and $J$ respectively.
Further we let $\mu^{(I)},\mu^{(J)}$ denote the respective marginal
distributions as defined in (\ref{marginal}). Then the condition for 
negative association may be written 
\[
\sum_{x_I,x_J} f(x_I)g(x_J) \mu(x_I,x_J) \leq \sum_{x_I,x_J} f(x_I)g(x_J)
\mu^{(I)}(x_I)\mu^{(J)}(x_J)
\]
or
\begin{equation}\label{NABoolean}
\sum_{x_I,x_J} f(x_I)g(x_J)\big(\mu(x_I,x_J)-\mu^{(I)}(x_I)\mu^{(J)}(x_J)\big)\leq 0.
\end{equation}
As in (\ref{NCPolynomial}), equation (\ref{NABoolean}) 
may be re-formulated as 
\begin{equation}\label{polynomialNA}
p_{f,g}(\mu_1,\dots,\mu_{2^n}) \leq 0,
\end{equation}
where $p_{f,g}(\mu_1,\dots,\mu_{2^n})$ is a polynomial in 
$\mu_1,\dots, \mu_{2^n}$ dependent on $f$ and $g$.

Equation (\ref{NCBoolean}) must hold for $1\leq i,j\leq n$, which 
is a finite number of constraints. 
Equation (\ref{NABoolean}) must hold for every non-decreasing $f$ and $g$, and 
disjoint index sets $I$ and $J$ -- an infinite number of constraints. Let
us however restrict our attention
to the set of all {\it strictly} NA distributions, i.e. those $\mu$ for which 
strict inequality holds in (\ref{NABoolean}):
\[
\C_{\mu}(f(x_I),g(x_J)) <0,
\]
for all monotone $f$ and $g$.
Multiplying (\ref{NABoolean}) 
by a constant, we may assume such $f$ and $g$ are uniformly bounded.
Note also that every function $f:\{0,1\}^n\to \R$ is a polynomial of bounded degree.
Consider the space of non-decreasing, uniformly bounded polynomials of degree at most $n$ (in the space
of continuous functions $C(\{0,1\}^n)$ equipped with the $\infty$-norm). 
Such spaces of polynomials of a finite degree compose a finite dimensional space. Moreover the supremum norm on this space is equivalent to the supremum norm on the coefficients, which is equivalent to the equicontinuity of the space. 
Thus by the Arzel\'a-Ascoli thoerem, this space is compact in the $\infty$-norm.
Combining this with the fact that the covariance operator is continuous,
there exists $\epsilon>0$ and finitely many $f_1,\dots, f_m$ and 
$g_1,\dots, g_m$ such that 
\[
\C_{\mu}(f(x_I),g(x_J)) <0
\]
$\forall$ $f,g$ non-decreasing, if 
\[
\C_{\mu}(f_i(x_I),g_i(x_J)) < -\epsilon
\]
$\forall$ $i,j=1,\dots, m$, or in the language of (\ref{polynomialNA}),
\begin{equation}\label{NAPolynomial2}
p_{f_i,g_i}(\mu_1,\dots,\mu_{2^n}) < -\epsilon, \;\;\;\;\; i=1,\dots,m.
\end{equation}

From the viewpoint of (\ref{NCBoolean}) and (\ref{NAPolynomial2}), 
the conditions of strict negative correlation and strict negative association on the 
Boolean cube are continuous in the parameters $\mu_1,\dots, \mu_{2^n}$
of a given distribution $\mu$. That is, the condition will still be satisfied
under small perturbations of $\mu_1,\dots, \mu_{2^n}$. 

Of course, the space of probability measures on the Boolean cube is a finite dimensional
space, and all Hausdorff vector topologies on a finite dimensional space are equivalent. Thus
one may define basic open sets by (\ref{openset}) or by (\ref{TVdiscrete}).
Or, equivalently, one may choose the Euclidean topology induced by the coordinate
system $\mu=(\mu_1,\dots,\mu_{2^n})$. 
As the conditions
defining both the class of NC and NA distributions are continuous in
this Euclidean topology, we moreover obtain,
\begin{thm}
Let $\m_{NC}$ and $\m_{NA}$ denote the spaces of NC and NA
distributions on the Boolean cube. We have
\[
\partial\m_{NC}\subset \{ \mu \in \m_{NC} : \mu^{(i,j)}(1,1) = \mu^{(i)}(1)\mu^{(j)}(1)
\text{ for some } i,j\}
\]
and
\[
\partial\m_{NA} \subset \{ \mu \in \m_{NA} : \exists f,g \text{ non-constant, non-decreasing, } \C_{\mu}(f(x_I),g(x_J)) = 0 \},
\]
where $I$ and $J$ are disjoint subsets of $\{1,\dots, n\}$.
Moreover, the interior of $\m_{NC}$ and the interior of $\m_{NA}$ are non-empty.
\end{thm}

\subsection{Convexity and Connectedness}
We further our study of the topological properties of the spaces of negatively associated and negatively correlated distributions by considering properties of convexity and connectedness. We consider such questions on both the Boolean cube and on all of $\R^n$.


\subsubsection{Convexity Properties of the Space of Negatively Associated Distributions}


\begin{thm}\label{convexNA}
The space of negatively associated
distributions is not convex on $\R^n$. 
\end{thm}

\begin{Proof}
We consider strictly negatively associated distributions $\mu$ and 
$\nu$, which exist by Lemma \ref{strictNAexist}. We show that there exist increasing functions $f,g$ and a $\lambda\in (0,1)$ for which the condition for negative association fails under the measure $\lambda\mu +(1-\lambda)\nu$. 

We begin with a general algebraic manipulation. Given increasing $f$ and $g$ defined on disjoint index sets,
there exist $\epsilon_1$ and $\epsilon_2$ such that 
\[
\int fg \, d\mu = \int f \, d\mu \int g \, d\mu -\epsilon_1
\]
and 
\[
\int fg \, d\nu = \int f \, d\nu \int g \, d\nu -\epsilon_2.
\]
Setting $A:= \int f\, d\mu$, 
$B:= \int f \, d\nu$, $C:= \int g \, d\mu$,
and $D := \int g\, d\nu$,
it follows that 
\begin{eqnarray*}
\int fg \, d(\lambda\mu+(1-\lambda)\nu) &=& \lambda\int fg\,d\mu +(1-\lambda)\int fg \, d\nu \\
&=& \lambda\int f \,d\mu \int g\,d\mu +(1-\lambda)\int f\, d\nu \int g \, d\nu 
-(\lambda\epsilon_1+(1-\lambda)\epsilon_2)\\
&=& \lambda AC+(1-\lambda)BD-(\lambda\epsilon_1+(1-\lambda)\epsilon_2).
\end{eqnarray*}
Further we have
\begin{eqnarray*}
\int f \, d(\lambda\mu+(1-\lambda)\nu) \int g \, d(\lambda\mu+(1-\lambda)\nu) &=&
(\lambda A +(1-\lambda)B)
(\lambda C +(1-\lambda)D)\\
&=&
\lambda^2 AC + \lambda(1-\lambda)AD
+\lambda(1-\lambda)BC + (1-\lambda)^2BD.
\end{eqnarray*}
The condition for convexity therefore becomes,
\[
\lambda AC +(1-\lambda)BD
-(\lambda\epsilon_1+(1-\lambda)\epsilon_2)
\leq \lambda^2 AC + \lambda(1-\lambda)AD
+\lambda(1-\lambda)BC + (1-\lambda)^2BD
\]
for $0\leq \lambda \leq 1.$ Simplifying, we obtain
\[
\lambda^2(A-B)(C-D) -\lambda(A-B)(C-D) \geq 
-(\lambda\epsilon_1+(1-\lambda)\epsilon_2).
\]
Upon setting $\tilde{C} = (A-B)(C-D)$ this becomes
\begin{equation}\label{convexity}
\tilde{C}\lambda^2 - \tilde{C} \lambda\geq 
-(\lambda\epsilon_1+(1-\lambda)\epsilon_2).
\end{equation}


Thus we must show that (\ref{convexity}) fails for certain increasing $f,g$ and $\lambda\in (0,1)$. If $\tilde{C}<0$, then the quadratic $\tilde{C}\lambda^2-\tilde{C}\lambda=\tilde{C}\lambda(\lambda-1)$ is non-negative for
all $0\leq \lambda \leq 1$, thus satisfying (\ref{convexity}). However, we claim that if 
$\tilde{C}>0$, then (\ref{convexity}) will not hold for certain $\lambda\in (0,1)$, as is shown below (this would complete the proof of non-convexity). 

To this end, we will need to guarantee $\tilde{C}>0$; this can be accomplished by choosing $\mu$ and $\nu$ in such a way that $A=\int f d\mu> \int f d\nu =B$, and $C=\int g d\mu >\int g d\nu =D$.
Given real numbers $p_1,\dots, p_n$ we can translate the mean of a probability measure $\mu$ in each variable by $p_i$ without changing the covariance structure of $\mu$. Specifically, map $\mu\mapsto \mu\circ T^{-1}$, where the transformation $T:\R^n\to \R^n$ is defined by 
\[
T(x_1,\dots, x_n) = (x_i+p_i)_{i=1,\dots,n} =(y_i)_{i=1,\dots,n}.
\]
By the change of variables formula, for any integrable $f$,
\[
\int_{\R^n} f(y) \, d(\mu\circ T^{-1})(y) = 
\int_{\R_n} f(Tx) \,d\mu(x) 
\]
In particular,
\[
\int_{\R^n} y_i \, d(\mu\circ T^{-1})(y) = \int_{\R^n}x_i\, d\mu(x) + p_i.
\]
Thus $\C_{\mu\circ T^{-1}}(y_i,y_j) = \C_{\mu}(x_i,x_j)$. What's more,
since $T$ is mere translation, it preserves the product ordering on $\R^n$. That is $x\geq y$ if and only if $Tx\geq Ty$. Therefore a function $f$ is non-decreasing on $\R^n$ if and only if $f\circ T$ is non-decreasing on $\R^n$, whence $\mu$ is negatively associated if and only if $\mu \circ T^{-1}$ is negatively associated.  
We may thus assume that $A=\int f d\mu> \int f d\nu =B$, and $C=\int g d\mu >\int g d\nu =D$. Thus $\tilde{C}>0$.

In this case the quadratic  $\tilde{C}\lambda^2-\tilde{C}\lambda$ is 
bounded above by 0 for all $0\leq \lambda \leq 1$, and its minimum value is attained at $\lambda=1/2$. If convexity is to hold, then (\ref{convexity}) must be valid when $\lambda=1/2$. Setting $\lambda=1/2$ in (\ref{convexity}) we obtain
\[
\frac{\tilde{C}}{4}\leq \frac{1}{2}(\epsilon_1+\epsilon_2).
\]
This is evidently false if $\epsilon_1$ and $\epsilon_2$ can be made arbitrarily small independent of $\tilde{C}$. 
We have constructed $\tilde{C}$ via a translation operator $T$ which is independent of $\epsilon_1$ and $\epsilon_2$, thus it suffices to show that these $\epsilon$-quantities can be made arbitrarily small. 

Following the proof of Lemma 1, we describe strictly negatively associated distributions $\mu$ and $\nu$. We consider once again the inequality 
\[
\frac{\tilde{C}}{4} \leq \frac{1}{2}(\epsilon_1+\epsilon_2)
\]
relative to the context at hand. It suffices to verify that for certain $f,g$ and measures $\mu, \nu$, said inequality is invalid. As in lemma 1, we may assume that the quantity $E_{\mu}fg =0$, and $E_{\nu}fg = 0$. It then suffices to verify $\epsilon_1$ and $\epsilon_2$ as being small. These quantities describe the level of negative association of the respective measures, for given functions $f$ and $g$. As said measures are supported on the standard basis vectors, $\epsilon_1$ and $\epsilon_2$ will take the form
\[
0-\left(\sum_k \beta_k f(\alpha_k)\right)\left( \sum_j \beta_j g(\alpha_j)\right),
\]
where 0 denotes the value of $E fg$, and $Ef = \sum_i\beta_i f(\alpha_i)$, relative to the measure $\mu = \sum_i\beta_i \delta_{\alpha_i}$. The analogue for $g$ and $\nu$. As $f$ and $g$ are arbitrary, it is clear that we may describe $\tilde{C}$ as positive. We need only guarantee that $f$ and $g$ are non-decreasing, while retaining that $\tilde{C}>\epsilon$ for some fixed $\epsilon>0$. In doing so, we let $Ef$ and $Eg$ approach 0 from above, thereby attaining arbitrarily small values for $\epsilon_1$ and $\epsilon_2$, and contradicting the inequality 
\[
\frac{\tilde{C}}{4} \leq \frac{1}{2}(\epsilon_1+\epsilon_2),
\]
as required. This shows that (\ref{convexity}) fails for $\lambda=1/2$, whereby convexity is violated for the corresponding distribution.
Thus it is shown that the space of negatively associated distributions on the Boolean cube is non-convex.

Thus it is shown that the space of negatively associated distributions on $\R^n$ is non-convex.
\end{Proof}

\subsubsection{Non-Convexity of $\m_{NC}(\R^n)$}
Define
the sets 
\begin{equation}\label{convexset}
E_{p_1,\dots, p_n} := \{\mu\in \m_{NC}(\R^n) :
\E_{\mu} x_i = p_i, i=1,\dots, n\}.
\end{equation}

\begin{cor}\label{NCCONVEX}
The space of negatively correlated distributions is not convex on $\R^n$. However, for any fixed $p_1,\dots, p_n\in \R$, the collection of measures
\[
E_{p_1,\dots,p_n} = \{\mu\in \m_{NC}(\R^n) : \E_{\mu} x_i = p_i,\; i=1,\dots,n\}
\]
is convex. 
\end{cor}

\begin{Proof}
The proof of non-convexity follows as in the proof of Theorem \ref{convexNA}. Specifically, 
given strictly negatively correlated distributions $\mu$ and $\nu$, fix $i,j$ and note that 
\[
\int x_ix_j \, d\mu = \int x_i \, d\mu \int x_j \, d\mu -\epsilon_1
\]
for some $\epsilon_1>0$, and 
\[
\int x_ix_j \, d\nu = \int x_i \, d\nu \int x_j \, d\nu -\epsilon_2
\]
for some $\epsilon_2>0$.
Set $A= \int x_i\, d\mu$, $B= \int x_i \, d\nu$, $C= \int x_j \, d\mu$,
and $D = \int x_j\, d\nu$. Then if convexity is to hold, we once again must have
\[
\lambda AC +(1-\lambda)BD
-(\lambda\epsilon_1+(1-\lambda)\epsilon_2)
\leq \lambda^2 AC + \lambda(1-\lambda)AD
+\lambda(1-\lambda)BC + (1-\lambda)^2BD
\]
for $0\leq \lambda \leq 1.$ Now with $\tilde{C}=(A-B)(C-D)$ this simplifies to 
\begin{equation}\label{Ceq}
\tilde{C}\lambda^2-\tilde{C}\lambda \geq -(\lambda \epsilon_1+(1-\lambda)\epsilon_2).
\end{equation}

If we can show that there exist NC distributions $\mu$ and $\nu$ such that $\tilde{C}>0$, then upon setting $\lambda = 1/2$ in (\ref{Ceq}), we will arrive at the condition 
\[
\frac{\tilde{C}}{4} \leq \frac{1}{2}(\epsilon_1+\epsilon_2)
\]
which will fail for small enough $\epsilon_1,\epsilon_2$, if we can make $\tilde{C}$ large enough independent of $\epsilon_1,\epsilon_2$. 

Thus we must show that there exist NC $\mu$ and $\nu$ such that $\tilde{C}>0$, and such that $\epsilon_1$ and $\epsilon_2$ are sufficiently small. That is 
$\int x_i\, d\mu > \int x_i\, d\nu$ and $\int x_j\, d\mu>\int x_j\, d\nu$, for $i\neq j$. Given real numbers $p_1,\dots, p_n$ we can translate the mean of a probability measure $\mu$ in each variable by $p_i$ without changing the covariance structure of $\mu$. Specifically, map $\mu\mapsto \mu\circ T^{-1}$, where the transformation $T:\R^n\to \R^n$ is defined by 
\[
T(x_1,\dots, x_n) = (x_i+p_i)_{i=1,\dots,n} =(y_i)_{i=1,\dots,n}.
\]
By the change of variables formula, for any integrable $f$,
\[
\int_{\R^n} f(y) \, d(\mu\circ T^{-1})(y) = 
\int_{\R_n} f(Tx) \,d\mu(x) 
\]
In particular,
\[
\int_{\R^n} y_i \, d(\mu\circ T^{-1})(y) = \int_{\R^n}x_i\, d\mu(x) + p_i.
\]
Thus $\C_{\mu\circ T^{-1}}(y_i,y_j) = \C_{\mu}(x_i,x_j)$. Thus given NC distributions $\mu$ and $\nu$, we may always translate $\mu$ until its mean values in each coordinate, i.e. $\int x_i \, d\mu(x)$, dominate the mean values of $\nu$ in each coordinate. This does not change the covariance structure of $\mu$, and therefore preserves negative correlation, and in particular $\epsilon_1$ and $\epsilon_2$.

We now prove that the collection 
\[
E_{p_1,\dots,p_n} = \{\mu\in \m_{NC}(\R^n) : \E_{\mu} x_i = p_i,\; i=1,\dots,n\}
\]
is convex for each fixed $p_1,\dots, p_n \in \R$. This follows from equation (\ref{Ceq}). Indeed if $\tilde{C}=0$ then certainly (\ref{Ceq}) holds for each $0\leq \lambda\leq 1$, and therefore the collection of NC distributions which satisfy $\tilde{C}=0$ for each $1\leq i<j\leq n$ will be convex. We have $\tilde{C}=0$ whenever $\int x_i\, d\mu = \int x_i \, d\nu$, and therefore if $\int x_i\, d\mu = \int x_i \, d\nu=p_i$ for each $i=1,\dots, n$, then their convex combination $\lambda \mu + (1-\lambda)\nu$ will be NC for each $0\leq \lambda\leq 1$. The result follows. 
\end{Proof}

\subsubsection{Non-Convexity of $\m_{NC}(I_n)$}
\begin{lemma}\label{NCmarginal}
For any $0<\epsilon<1$ and $1\leq i\leq n$, there exists a negatively correlated  distribution on the Boolean cube $I_n=\{0,1\}^n$ satisfying $\mu^{(i)}(1)=\epsilon$. In fact, given $0<\epsilon_i<1$, $i=1,\dots, n$, satisfying $\sum_i \epsilon_i=1$, there exist a negatively correlated distribution on the Boolean cube satisfying $\mu^{(i)}(1)= \epsilon_i$, $i=1,\dots, n$.
\end{lemma}

\begin{Proof}
Given $1\leq i\leq n$, define a distribution as follows:
\[
\mu_{i,\epsilon} := \epsilon\delta_{\alpha_i} +(1-\epsilon)\delta_{\alpha_j},
\]
where $i\neq j$ and ${ \alpha_i} = (0,\dots,1,\dots,0)$ is the vector with a single 1 in the $i$th component (likewise for $\alpha_j$). We see that $\mu^{(i)}_{i,\epsilon}(1) =\epsilon$. What's more, we have
\[
{\bf E}_{\mu_{i,\epsilon}} x_{k}x_{\ell} =\mu^{(k,\ell)}_{i,\epsilon}(1,1)= 0
\]
for all $1\leq k,\ell \leq n$, and thus
\[
{\bf E}_{\mu_{i,\epsilon}} x_{k}x_{\ell} - {\bf E}_{\mu_{i,\epsilon}} x_k{\bf E}_{\mu_{i,\epsilon}} x_{\ell} \leq 0.
\]
The first result follows.

Define a measure $\mu$ as the convex combination of point masses centered at each $\alpha_i =(0,\dots, 1,\dots, 0)$, $i=1,\dots, n$ (where $\alpha_i$ has exactly one 1 in the $i$th coordinate, as above):
\[
\mu : = \sum_ i \epsilon_i \delta_{\alpha_i}.
\]
Evidently, for each $j=1,\dots, n$, $\mu^{(j)}(1) =\epsilon_j$, and for each $i\neq j$ $\mu^{(i,j)}(1,1) = 0$. Thus $\mu$ is strictly negatively correlated:
\[
\mu^{(i,j)}(1,1) - \mu^{(i)}(1)\mu^{(j)}(1) = -\epsilon_i\epsilon_j.
\]
Further $\mu$ satisfies the requirements on the one-dimensional marginals. The result follows.
\end{Proof}

\begin{cor}
The space of negatively correlated distributions on the Boolean cube $I_n$ is non-convex. However, for any fixed $p_1,\dots, p_n \in \R$, the collection of measures
\[
E_{p_1,\dots, p_n} = \{\mu\in \m_{NC}(I_n) : \mu^{(i)}(1) = p_i, \; i=1,\dots, n\}
\]
is convex.
\end{cor}

\begin{Proof}
We begin with strictly negatively correlated distributions $\mu$ and 
$\nu$, which exist by Lemma \ref{strictNAexist}. We derive conditions under which the convex combination $\lambda\mu +(1-\lambda)\nu$ fails to be negatively correlated. We then produce strictly negatively correlated distributions whose convex combination fails to satisfy the above-mentioned condition.

Note that $\E_{\mu} x_i =\mu^{(i)}(1)$ and $\E_{\mu}x_ix_j =\mu^{(i,j)}(1,1)$ on the Boolean cube. Thus given strictly negatively correlated distributions $\mu$ and $\nu$, fix $i,j$ and note that 
\[
\mu^{(i,j)}(1,1) = \mu^{(i)}(1)\mu^{(j)}(1) -\epsilon_1
\]
for some $\epsilon_1>0$, and 
\[
\nu^{(i,j)}(1,1) = \nu^{(i)}(1)\nu^{(j)}(1) -\epsilon_2
\]
for some $\epsilon_2>0$.
Set $A= \mu^{(i)}(1)$, $B= \nu^{(i)}(1)$, $C= \mu^{(j)}(1)$,
and $D = \nu^{(j)}(1)$. Then if convexity is 
to hold, we once again must have (see proof of Theorem \ref{convexNA} above, where the roles of $f$ and $g$ are played by $x_i$ and $x_j$)
\[
\lambda AC +(1-\lambda)BD
-(\lambda\epsilon_1+(1-\lambda)\epsilon_2)
\leq \lambda^2 AC + \lambda(1-\lambda)AD
+\lambda(1-\lambda)BC + (1-\lambda)^2BD
\]
for $0\leq \lambda \leq 1.$ 
Now with $\tilde{C}=(A-B)(C-D)$ this simplifies to 
\[
\tilde{C}\lambda^2-\tilde{C}\lambda \geq -(\lambda \epsilon_1+(1-\lambda)\epsilon_2).
\]
As this holds for $\tilde{C}=0$, the second statement of the Corollary holds. 

Now set $\lambda=1/2$. We arrive at the condition 
\[
\frac{\tilde{C}}{4} \leq \frac{1}{2}(\epsilon_1+\epsilon_2).
\]
Our strategy is as follows. We introduce measures $\mu$ and $\nu$ as convex combinations of point masses at the standard basis vectors. We demonstrate that the inequality
\[
\frac{\tilde{C}}{4} \leq \frac{1}{2}(\epsilon_1+\epsilon_2)
\]
is valid. We then define a perturbation of this measure under which said inequality is violated, whereby we obtain a convex combination of measures which fail to be negatively correlated.

According to Lemma \ref{NCmarginal}, there exist strictly negatively correlated distributions with $\tilde{C}>0$. Specifically, define 
\[
\mu = \sum_k \beta_k \delta_{\alpha_k},
\]
and further
\[
\nu = \sum_k \beta_k' \delta_{\alpha_k}
\]
where $\alpha_k = (0,\dots, 1,\dots, 0)$ has exactly one 1 in the $k$th component, and $\sum _k \beta_k =1$, $\sum_k \beta_k'=1$. 
Then 
\[
\tilde{C} = (\beta_i -\beta_i')(\beta_j - \beta_j')
\]
where we have used that  $\tilde{C}=(A-B)(C-D)$ and $A= \mu^{(i)}(1)$, $B= \nu^{(i)}(1)$, $C= \mu^{(j)}(1)$,
and 
\[
\epsilon_1 = \beta_i\beta_j,
\;\;\;\epsilon_2 =\beta_i'\beta_j'.
\]
Thus the condition 
\[
\frac{\tilde{C}}{4} \leq \frac{1}{2}(\epsilon_1+\epsilon_2)
\]
becomes 
\[
\frac{(\beta_i -\beta_i')(\beta_j - \beta_j')}{4} \leq \frac{1}{2}(\beta_i\beta_j+\beta_i'\beta_j').
\]
This inequality reduces to 
\[
\beta_i\beta_j +\beta_i'\beta_j' +\beta_i\beta_j' + \beta_i'\beta_j\geq 0,
\]
which holds for all non-negative reals.  

Since strict negative correlation is continuous in the Euclidean parameters of the distribution, we may perturb the above convex measure as follows. Define $\mu$ as above, but perturbed with a small additional weight given to $(1,1,0,\dots, 0)$: $\mu(1,1,0,\dots, 0) = \epsilon$. Here $\epsilon>0$ is small enough so that $\mu$ is still negatively correlated.  
Compensating this increased weight by a total decrease in the other decoupled positive weights totaling $\epsilon$ to keep normalization will still not affect negative correlation if $\epsilon$ is sufficiently small.
Then $\mu^{(1,2)}(1,1) = \epsilon$ and $\mu^{(i)}(1) = \beta_i+\epsilon$ for $i=1,2$ . If we do the same for $\nu$, with the same $\epsilon$ perturbation, then note that $\tilde{C}$ does not change. Indeed, for $i=1, j=2$ 
\begin{eqnarray*}
    \tilde{C} &=& (A-B)(C-D) \\
    &=& (\beta_1+\epsilon-\beta_1'-\epsilon)(\beta_2+\epsilon -\beta_2'-\epsilon)\\
    &=&(\beta_1-\beta_1')(\beta_2-\beta_2').
    \end{eqnarray*}
We may assume that $\tilde{C}>0$, by virtue of choosing $\beta_1>\beta_1'$ and $\beta_2>\beta_2'$.
On the other hand, for $\epsilon>0$ small enough
\[
 -\epsilon_1 \equiv \mu^{(1,2)}(1,1) -\mu^{(1)}(1)\mu^{(2)}(1) <0
\]
as the quantity 
\[
\mu^{(1,2)}(1,1) -\mu^{(1)}(1)\mu^{(2)}(1) =\epsilon-(\beta_1+\epsilon)(\beta_2+\epsilon)
\]
is continuous in $\epsilon$, and approaches $-\beta_1\beta_2$ from above as $\epsilon\to 0$. For $\beta_1$ and $\beta_2$ small enough, the quantity 
\[
\epsilon-(\beta_1+\epsilon)(\beta_2+\epsilon)
\]
will be positive; letting $\epsilon$ tend to 0 shows that this quantity will pass through zero, towards $-\beta_1\beta_2$. This shows that $\epsilon_1$ will move through 0, and thus the inequality 
\[
\frac{\tilde{C}}{4} \leq \frac{1}{2}(\epsilon_1+\epsilon_2)
\]
will be violated, as consideration of $\epsilon_2$ will be analogous. 

It follows that the right hand side of the inequality 
\[
\frac{\tilde{C}}{4} \leq \frac{1}{2}(\epsilon_1+\epsilon_2)
\]
will decrease. And thus we can, by continuity, decrease the right hand side until the inequality is violated. Thus it is shown that the space of negatively correlated distributions on the Boolean cube is non-convex.
\end{Proof}

\subsubsection{Non-Convexity of $\m_{NA}(I_n)$}     

\begin{lemma}
 For any $1\leq i\leq n$ and $0<\epsilon<1$, there exists a negatively associated distribution $\mu$ on the Boolean cube $I_n=\{0,1\}^n$ satisfying $\mu^{(i)}(1)=\epsilon$.
\end{lemma}

\begin{Proof}
As in the proof of Lemma \ref{strictNAexist}, any distribution supported on \[
I_{n,1}=\left\{(x_1,\dots,x_n)\in I_n: \sum_j x_j =1\right\} 
\]
is strictly negatively associated. Such a measure will be of the form
\[
\mu=\sum_{{\bf x}_k \in I_{n,1}} \alpha_k \delta_{{\bf x}_k}
\]
where $\sum_k \alpha_k=1$. Given $1\leq i\leq n$ and $0<\epsilon<1$, in setting $\alpha_i =\epsilon$, and determining the remaining coefficients by the condition $\sum_k \alpha_k =1$, we see that $\mu^{(i)}(1)=\epsilon$.
\end{Proof}
                   
\begin{cor}
The space of negatively associated distributions on the Boolean cube $I_n$ is non-convex.
\end{cor}   
    
\begin{Proof}
Following once again the proof of Theorem \ref{convexNA}, we set  $A:= \int_{I_n} f\, d\mu$, 
$B:= \int_{I_n} f \, d\nu$, $C:= \int_{I_n} g \, d\mu$,
and $D := \int_{I_n} g\, d\nu$, and obtain the same inequality dictating convexity:
\[
\lambda^2(A-B)(C-D) -\lambda(A-B)(C-D) \geq 
-(\lambda\epsilon_1+(1-\lambda)\epsilon_2)
\]
for all $0\leq \lambda\leq 1$. Recall once again that $\epsilon_1$ and $\epsilon_2$ are defined by the initial assumption: Given increasing $f$ and $g$ defined on disjoint index sets,
there exist $\epsilon_1$ and $\epsilon_2$ such that 
\[
\int fg \, d\mu = \int f \, d\mu \int g \, d\mu -\epsilon_1
\]
and 
\[
\int fg \, d\nu = \int f \, d\nu \int g \, d\nu -\epsilon_2.
\]

Upon setting $\tilde{C} = (A-B)(C-D)$ this becomes
\[
\tilde{C}\lambda^2 - \tilde{C} \lambda\geq 
-(\lambda\epsilon_1+(1-\lambda)\epsilon_2).
\]
This condition is satisfied whenever $\tilde{C}=0$ or $\tilde{C}<0$, and fails when $\tilde{C}>0$ for small enough $\epsilon>0$.

Following the proof of Lemma 1, we describe strictly negatively associated distributions $\mu$ and $\nu$. We consider once again the inequality 
\[
\frac{\tilde{C}}{4} \leq \frac{1}{2}(\epsilon_1+\epsilon_2)
\]
relative to the context at hand. It suffices to verify that for certain $f,g$ and measures $\mu, \nu$, said inequality is invalid. As in lemma 1, we may assume that the quantity $E_{\mu}fg =0$, and $E_{\nu}fg = 0$. It then suffices, as in the proof of the previous proposition, to verify $\epsilon_1$ and $\epsilon_2$ as being small. These quantities describe the level of negative association of the respective measures, for given functions $f$ and $g$. As said measures are supported on the standard basis vectors, $\epsilon_1$ and $\epsilon_2$ will take the form
\[
0-\left(\sum_k \beta_k f(\alpha_k)\right)\left( \sum_j \beta_j g(\alpha_j)\right),
\]
where 0 denotes the value of $E fg$, and $Ef = \sum_i\beta_i f(\alpha_i)$, relative to the measure $\mu = \sum_i\beta_i \delta_{\alpha_i}$. The analogue for $g$ and $\nu$. As $f$ and $g$ are arbitrary, it is clear that we may describe $\tilde{C}$ as positive. We need only guarantee that $f$ and $g$ are non-decreasing, while retaining that $\tilde{C}>\epsilon$ for some fixed $\epsilon>0$. In doing so, we let $Ef$ and $Eg$ approach 0 from above, thereby attaining arbitrarily small values for $\epsilon_1$ and $\epsilon_2$, and contradicting the inequality 
\[
\frac{\tilde{C}}{4} \leq \frac{1}{2}(\epsilon_1+\epsilon_2),
\]
as required. This shows that (\ref{convexity}) fails for $\lambda=1/2$, whereby convexity is violated for the corresponding distribution.
Thus it is shown that the space of negatively associated distributions on the Boolean cube is non-convex.
\end{Proof}  


\subsubsection{Connectedness Properties of the Spaces of Negatively Correlated and Negatively Associated Distributions}        
\begin{thm}
The space of negatively correlated and the space of negatively associated distributions on the Boolean cube, and on $\R^n$, are path connected in the weak topology.

\end{thm} 
\begin{Proof}
For any negatively associated measure $\mu$, consider the family of measures $\mu_t$ such that for any set $A$, $\mu_t(A) = \mu(A/t)$; here $0< t\leq 1$. 
For any negatively associated measure $\mu$, the corresponding $\mu_t$ is negatively associated. 
For $t=0$ we define $\mu_0$ to be the point mass at the origin. We have defined $A/t$ to be the set of all points in $A$ divided by the constant $t$. As we scale $t$ from 1 to 0, this effectively concentrates the measure $\mu$ through this scaling into a point mass at the origin, while preserving negative association in the process. Indeed, assuming $A$ is a ball away from the origin, the mass of the set $A/t$ approaches zero as $t\to 0$, as the distance of the set $A/t$ from the origin approaches infinity as $t\to 0$. It follows that $\mu(A/t)$ converges weakly to the point mass at the origin. This provides a path connecting any two negatively associated distributions to the point mass at 0, proving that the family is path connected in the weak topology. 
\end{Proof}

\bibliographystyle{plain}
\bibliography{Thesis_bib}

\end{document}